\documentclass[11pt]{article}
\usepackage{amsmath}
\usepackage{amssymb}
\usepackage{amsfonts}

\newtheorem{theo}{Theorem}

\newtheorem{lemma}{Lemma}

\def\bydef{\stackrel{\mathrm{def}}{=}}%-----------------by definition
\newcommand{\dept}[1]{\par\pagebreak[2]\vspace{5mm}%
{{\raggedright \bf\hspace{0pt}#1\par}}%
\nopagebreak\vspace{2mm}}%-----------------------------------paragraph
\def\N{\mathbb{N}}% - natural
\def\R{\mathbb{R}}% - real
\def\C{\mathbb{C}}% - complex
\def\I{\mathrm{I}}%---positive-negative-neutral index sets
\def\G{\mathrm{G}}%-----domain of good values of alpha
\def\H{\mathrm{H}}
\def\L{\mathrm{L}}

\def\F{{\cal F}}%----------general Fischer-type
\def\B{{\cal B}}%----------general Bergman-type
\title{Hypergeometric reproducing kernels and analytic continuation from a half-line}
\author{Dmitrii Karp}
\date{}
\begin{document}
\maketitle \thispagestyle{myheadings} \markright{ Journal of
Integral Transforms and Special Functions, vol.14, no.6, 2003, pp.
485-498}

\begin{center}
Laboratory of Function Theory,\\ Institute of Applied Mathematics,\\
Far Eastern Branch of the Russian Academy of Sciences,\\ 7 Radio
Street, Vladivostok 690043, Russian Federation\\
%Department of Mathematics and Modeling,\\ Far Eastern Academy of
%Economics and Management,\\ Okeanksii Prospekt 19, Vladivostok,\\
%Russian Federation\\
E-mail: garfy@mail.primorye.ru
\end{center}

\vspace{0.4in}
\begin{center}
\parbox{12cm}{
\small\textbf{Abstract.} Indefinite inner product spaces of entire
functions and functions analytic inside a disk are considered and
their completeness studied.  Spaces induced by the rotation invariant
reproducing kernels in the form of the generalized hypergeometric
function are completely characterized.  A particular space generated
by the modified Bessel function kernel is utilized to derive an
analytic continuation formula for functions on $\R^+$ using the best
approximation theorem of S. Saitoh. As a by-product several new area
integrals involving Bessel functions are explicitly evaluated. }
\end{center}



\dept{1. Bergman spaces with indefinite inner product} The Bergman
space $A^2_w$ with radial weight $w$ on the unit disk $\Delta$ is
defined as the set of functions analytic inside $\Delta$ equipped
with inner product
\begin{equation}\label{ClBergman}
[f,g]_w=\int\limits_{\Delta}f(z)\overline{g(z)}w(|z|)d\sigma(z),
\end{equation}
where $d\sigma$ is the planar Lebesgue measure,
$w:(0,1)\to\R^{+}\cup\{0\}$ and $xw\in L_1(0,1)$. It seems to be the
folklore result that this space is complete iff
\begin{equation}\label{WeightCond}
\int\limits_{\xi}^{1}\!\!w(x)dx>0,~~~~~\forall\xi\in(0,1).
\end{equation}
It is interesting to remark here that for non-radial non-negative
weight no necessary and sufficient conditions for completeness seem
to be known. Now we relax the assumption of non-negativeness of the
weight $w$ and consider the space $A^2_w$ on the disk $\Delta_R$ with
the (possibly infinite) radius $R$  determined by $w$ as follows
\begin{equation}\label{Rdef}
R\bydef\sup\left\{x:\int\limits_{x}^{\infty}|w(x)|dx>0\right\}.
\end{equation}
The space $A^2_w(\Delta_R)$ with sign-alternating weight $w$ is of
course no longer a Hilbert space. By standard argument (which carries
over with minor modifications to the  sign-alternating case) the
inner product (\ref{ClBergman}) in this space can be written as
\begin{equation}\label{InnerProd-Series}
\int\limits_{\Delta_R}\!\!f(z)\overline{g(z)}w(|z|)d\sigma=
2\pi\sum\limits_{k=0}^{\infty}f_k\overline{g}_k\int\limits_{0}^{R}r^{2k+1}w(r)dr=
\sum\limits_{k=0}^{\infty}c_k f_k\overline{g}_k,~~~ c_k\bydef
2\pi\int\limits_{0}^{R}r^{2k+1}w(r)dr.
\end{equation}
Define
\[
\I_+\bydef\{k:c_k>0\},~~\I_-\bydef\{k:c_k<0\},~~
\I_0\bydef\{k:c_k=0\}.
\]
The space $A^2_w(\Delta_R)$ allows the fundamental decomposition
\begin{equation}\label{Decomp}
A^2_w=A^2_{w+}\oplus A^2_{w-},~~~f=f_+ + f_-,
\end{equation}
with
\begin{equation}\label{Decomp1}
f_+(z)=\sum\limits_{k\in \I_+}f_kz^k,~~~f_-(z)=\sum\limits_{k\in
 \I_-}f_kz^k,
\end{equation}
and
\[
[f_+,f_+]=\sum\limits_{k\in\I_+}c_k|f_k|^2>0,~~[f_-,f_-]=\sum\limits_{k\in\I_-}c_k|f_k|^2<0,
~~[f_+,f_-]=0.
\]
Thus $A^2_w(\Delta_R)$ is a Krein space (i.e. complete indefinite
inner product space allowing decomposition (\ref{Decomp}), cf.
\cite{Alpay2,AI,Bognar}) if the spaces $A^2_{w+}$ and $A^2_{w-}$ are
both complete. Define
\begin{equation}\label{ak}
a_k\bydef 2\pi\int\limits_{0}^{R}x^{2k+1}|w(x)|dx.
\end{equation}
\begin{theo}\label{CompleteKrein}
The space $A^2_w$ is a (complete) Krein space if and only if the
sequence  $\{a_k/c_k\}_{k\in\I_{-}\cup\I_{+}}$ is bounded.
\end{theo}
The proof is straightforward and is given in \cite{Karp4}. One can
also find an example there of weight that violates the boundedness
condition from Theorem~\ref{CompleteKrein}.  The space $A^2_w$ admits
the reproducing kernel given by
\begin{equation}\label{ReprGeneral}
\varphi(z\overline{u})=\sum_{k\in\I_{-}\cup\I_{+}}\frac{1}{c_k}(z\overline{u})^k.
\end{equation}
For the general theory of reproducing kernels in Hilbert and Krein
spaces see \cite{Alpay2,Ar2,SaitohBook}. The condition of
Theorem~\ref{CompleteKrein} is obviously satisfied when there exists
$\gamma<R$ such that $w\geq 0$ or $w\leq 0$ when $x\in(\gamma,R)$. If
this requirement is met and all moments (\ref{InnerProd-Series}) are
finite we will call $w$ an {\em admissible weight} on $(0,R)$.
Moreover, in the case of admissible weight at least one of the sets
$\I_{+}$ or $\I_{-}$ is finite implying that at least one of the
spaces $A^2_{w+}$ or $A^2_{w-}$  is finite-dimensional and hence
$A^2_w(\Delta_R)$ becomes a Pontryagin space \cite{Alpay2,AI,Bognar}.
We call $A^2_w(\Delta_R)$ a positive Pontryagin space when
$\dim{A^2_{w-}}<\infty$ and a negative Pontryagin space when
$\dim{A^2_{w+}}<\infty$. The numbers $\dim{A^2_{w+}}$ and
$\dim{A^2_{w-}}$ are called positive and negative (Pontryagin)
indices of $A^2_w(\Delta_R)$, respectively.

\dept{2. Hypergeometric kernels} We now turn to the case when
$\varphi$ is the generalized hypergeometric function (cf.
\cite{Prud3}):
\begin{equation}\label{pFq-def}
\varphi(z\overline{u})= {_{p}F_{q}}\left(\left.\begin{array}{l}
a_1,\ldots, a_p
\\ b_1,\ldots, b_q
\end{array}\right|z\overline{u}\right)\bydef\sum\limits_{k=0}^{\infty}\frac{(a_1)_k\ldots(a_p)_k}
{(b_1)_k\ldots(b_q)_k}\frac{(z\overline{u})^k}{k!},
\end{equation}
where $(a)_0=1$, $(a)_k=a(a+1)\ldots(a+k-1)=\Gamma(a+k)/\Gamma(a)$ is
the Pochhammer symbol or shifted factorial. As shown below the weight
$w$ will be in this case expressed in terms of Meijer's G-function
defined by
\[
\G^{m,n}_{p,q}\!\!\left(\!z~\vline\begin{array}{l}a_1\ldots a_p
\\b_1\ldots b_q\end{array}\!\!\right)=
\G^{m,n}_{p,q}\left(\!z~\vline\begin{array}{l}(a_p)\\(b_q)\end{array}\!\!\right)\!\bydef
\]\[
\frac{1}{2\pi{i}}
\int\limits_{L_G}\!\!\frac{\Gamma(b_1\!+\!s)\dots\Gamma(b_m\!+\!s)\Gamma(1\!-\!a_1\!-\!s)\dots\Gamma(1\!-\!a_n\!-\!s)z^{-s}}
{\Gamma(a_{n+1}\!+\!s)\dots\Gamma(a_p\!+\!s)\Gamma(1\!-\!b_{m+1}\!-\!s)\dots\Gamma(1\!-\!b_q\!-\!s)}ds.
\]
Here $(m,n,p,q)$ is the order of G and $0\le m\le q,~ 0\le{n}\le p$.
The contour  $L_G$ is an infinite loop separating the poles of gamma
functions of the form $\Gamma(b_j+s)$, $j=1,2,\dots,m$ from those of
gamma functions of the form $\Gamma(1-a_j-s)$, $j=1,2,\dots,n$.
Details can be found in \cite{Lyakh,Prud3}.  We will only need
G-functions with $m=q$, $n=0$. On denoting
\begin{equation}\label{hs}
h(s)=\frac{\Gamma(b_1+s)\dots\Gamma(b_q+s)}{\Gamma(a_1+s)\dots\Gamma(a_p+s)},
\end{equation}
the residue theorem leads to the following expansion
\begin{equation}\label{Gexpress}
\G^{q,0}_{p,q}\left(\!x~\vline\begin{array}{l}(a_p)\\(b_q)\end{array}\!\!\right)\!=\!\!\!
\sum\limits_{\stackrel{j=\overline{1,q}}{\scriptscriptstyle{k=\overline{0,\infty}}}}
\frac{x^{b_j+k}}{(n_{jk}-1)!}\sum\limits_{i=0}^{n_{jk}-1}{\mathrm{C}}^i_{n_{jk}-1}\ln^{i}\frac{1}{x}
\left[(s+b_j+k)^{n_{jk}}h(s)\right]^{(n_{jk}-1-i)}_{|s=-b_j-k},
\end{equation}
where ${\mathrm{C}}^i_j$ are binomial coefficients, $n_{jk}$ is the
multiplicity of the pole of $h$ at $s_{j,k}=-b_j-k$ or $0$ if there
is no pole at the point $s_{j,k}$ and the convention $1/(-1)!=0$ is
used. The basic property of G-function is that it is the inverse
Mellin transform of $h$:
\begin{equation}\label{GMellin}
\int\limits_{0}^{\infty}\!\!x^{s-1}
\G^{q,0}_{p,q}\left(x~\vline\begin{array}{l}(a_p)\\(b_q)\end{array}\right)dx=
\frac{\Gamma(b_1+s)\ldots\Gamma(b_q+s)}{\Gamma(a_1+s)\ldots\Gamma(a_p+s)}.
\end{equation}

We will consider two cases $p<q$ and $p=q$ separately.
\begin{lemma}\label{G-asymp-lemma}
Let $p<q$, $b_i>-1$, $i=\overline{1,q}$, possibly except for $b_i$ of
the form $b_i=a_j+k$ for some $j\in\{1,2,\ldots,p\}$ and $k\in\N_0$;
$\theta>0$ and $a_i\in\R$. Then
$\G^{q,0}_{p,q}\left({\theta}x^2~\vline\begin{array}{l}(a_p)\\(b_q)\end{array}\right)$
is an admissible weight on $(0,\infty)$.
\end{lemma}
The proof hinges on the following asymptotic relations:
\begin{equation}\label{G-asymp-zero}
\G^{q,0}_{p,q}\left({\theta}x^2~\vline\begin{array}{l}(a_p)\\(b_q)\end{array}\right)\sim
x^{2b}\ln^{j-1}{x},~~~x\to{0},
\end{equation}
\begin{equation}\label{G-asymp}
\G^{q,0}_{p,q}\left(x~\vline\begin{array}{l}(a_p)\\(b_q)\end{array}\right)=
\frac{(2\pi)^{\frac{1}{2}(\mu-1)}}{\sqrt{\mu}}x^{(1-\alpha)/\mu}e^{-\mu
x^{1/\mu}}\left[1+O(x^{-1/\mu})\right],~~~x\to{\infty},
\end{equation}
where $\mu=q-p$,
$\alpha=\sum_{i=1}^{p}a_i-\sum_{i=1}^{q}b_i+\frac{1}{2}(q-p+1)$. Here
the symbol $\sim$ means that the ratio of the quantities on the left
and on the right is bounded from above an from below by positive
constants independent of $x$.  The first relation is a consequence of
(\ref{Gexpress}); the second is a particular case of the formula on
page 289 that follows (7.8) in \cite{Braaksma}.

Let $[x]$ denote the largest integer not exceeding $x$. Having made
these preparations we arrive at the following
\begin{theo}\label{KreinPlane}
Let  $(a_p)$, $(b_q)$ be real none of them being a non-positive
integer, $p\leq q$ are positive integers, $\theta>0$, $m$ is the
number of negative elements of the set $\{(a_p), (b_q)\}$, $d_i$ is
the negative element with  $i$-th smallest absolute value and
\[
\nu\bydef\sum\limits_{k=1}^{m}[d_k].
\]
Then the space $\F^{\,\theta}_{(a_p),(b_q)}$ induced by the
reproducing kernel {\em(\ref{pFq-def})} is the Pontryagin space of
entire functions with inner product
\begin{equation}\label{InnerPlaneSeries-intro}
[f,g]=\sum\limits_{k=0}^{\infty}\frac{(b_1)_k\ldots(b_q)_kk!}{(a_1)_k\ldots(a_p)_k\theta^k}
f_k\overline{g_k},
\end{equation}
where $f_k$, $g_k$ are Taylor's coefficients of $f$ and $g$,
respectively. It is a positive space with negative index
\begin{equation}\label{NumberNeg}
\sum\limits_{i=0}^{m-1}\sum\limits_{k=|[d_i]|+1}^{|[d_{i+1}]|}
\left(\sum\limits_{j=1}^{i}[d_j]+k(m-i)\right)\bmod{2},
\end{equation}
iff $\nu$ is even, and a negative space with positive index
\begin{equation}\label{NumberPos}
\sum\limits_{i=0}^{m-1}\sum\limits_{k=|[d_i]|+1}^{|[d_{i+1}]|}
\left[1-\left(\sum\limits_{j=1}^{i}[d_j]+k(m-i)\right)\bmod{2}\right]+1,
\end{equation}
iff $\nu$ is odd. Let
\[
l=\min\left\{0,[\min\limits_{b_i\ne a_j+k}b_i]\right\}.
\]
Then the inner product {\em(\ref{InnerPlaneSeries-intro})} may be
represented in the form
\[
[f,g]\!=\!\frac{\theta^{1-l}}{\pi}\Gamma\!\!\left[\!\begin{array}{l}
        (a_p) \\
        (b_q)
\end{array}
\!\right] \int\limits_{\C}\!\!\!f^{(l)}(z)\overline{g^{(l)}(z)}
\G^{q+2,0}_{p+1,q+2}\left(\theta|z|^2\left|\begin{array}{l}
        (a_p)+l-1,l \\
        (b_q)+l-1,0,0
\end{array}
\right.\!\!\!\right)\!d\sigma+
\]
\begin{equation}\label{InnerPlane1}
+\sum\limits_{k=0}^{l-1}\frac{(b_1)_k\ldots(b_q)_k}{(a_1)_k\ldots(a_p)_k
\theta^k k!}f^{(k)}(0)\overline{g^{(k)}(0)}.
\end{equation}
In particular if $l=0$
\begin{equation}\label{InnerPlane}
[f,g]=\frac{\theta}{\pi}\Gamma\!\!\left[\begin{array}{l}(a_p)\\(b_q)\end{array}\right]
\!\int\limits_{\C}\!\!f(z)\overline{g(z)}
\G^{q+1,0}_{p,q+1}\left(\!\theta|z|^2~\vline\begin{array}{l}(a_p)-1\\0,(b_q)-1\end{array}\!\!\right)
\!d\sigma.
\end{equation}
\end{theo}
\textbf{Proof.}  The formulae (\ref{NumberNeg}) and (\ref{NumberPos})
are derived by thorough counting of positive and negative
coefficients in (\ref{pFq-def}). Expressions (\ref{InnerPlane1}) and
(\ref{InnerPlane}) are obtained by comparing (\ref{InnerProd-Series})
and (\ref{ReprGeneral}) with (\ref{pFq-def}), (\ref{GMellin}) and
(\ref{InnerPlaneSeries-intro}). Existence of integral in
(\ref{InnerPlane1}) and (\ref{InnerPlane}) is guaranteed by
Lemma~\ref{G-asymp-lemma}. $\square$

When  $p=q=0$, the space $\F^{\,\theta}_{(a_p),(b_q)}$ becomes the
classical Bargmann-Fock space \cite{Barg1,Barg2}. When $p=q=a=1$,
$b>0$ this space reduces to the generalized  Bargmann-Fock space
considered in \cite{Burbea1,Burbea2}. Many more particular cases are
given in \cite{Karp4}.

For formulating a similar result in the case $p=q$ the following
lemma plays the key role.

\begin{lemma}\label{G-OutsideDisk}
If $x>1$ and
$\Re\left(\sum_{i=0}^{p}b_i\right)<\Re\left(\sum_{i=0}^{p}a_i\right)$,
then
\begin{equation}
\G^{p,0}_{p,p}\left(x\left|\begin{array}{l}
        \left(a_p\right)  \\
        \left(b_p\right)
\end{array}
\right.\right)=0.
\end{equation}
\end{lemma}

For a proof the reader is again referred to \cite{Karp4}.  Since
$\G^{p,0}_{p,p}(x)$ has a (finite or infinite) limit as $x\to{1}$, it
is an admissible weight on $(0,1)$. Thus we get
\begin{theo}\label{KreinDisk}
Let $(a_p)$, $(b_{p-1})$, $\theta>0$, $m$, $d_i$ and $\nu$ have the
same meaning as in Theorem~\ref{KreinPlane}. Then the space
$\B^{\,\theta}_{(a_p),(b_{p-1})}$ induced by the reproducing kernel
\begin{equation}\label{RKDisk}
K(z, \overline{u})={_{p}F_{p-1}}\left(\left.\begin{array}{l}
        a_1,\ldots, a_p  \\
        b_1,\ldots, b_{p-1}
\end{array}\right|\frac{1}{\theta}z\overline{u}\right)
\end{equation}
is a Pontryagin space comprising functions holomorphic inside
$\Delta_{\sqrt{\theta}}$ with inner product
\begin{equation}\label{InnerDiskSeries}
[f,g]=\sum\limits_{k=0}^{\infty}\frac{(b_1)_k\ldots(b_{p-1})_kk!\theta^k}
{(a_1)_k\ldots(a_p)_k} f_k\overline{g_k},
\end{equation}
where $f_k$, $g_k$ are Taylor's coefficients of $f$ and $g$,
respectively. This space is positive and has negative index given by
{\em(\ref{NumberNeg})} iff $\nu$ is even and is negative with
positive index given by {\em(\ref{NumberPos})} iff $\nu$ is odd.
Define
\[
l_1=\min\left\{0,[\min\limits_{b_i\ne a_j+k}b_i]\right\},
\]
and
\[
l_2=\min\left\{k\in\N_0:\sum\limits_{i=0}^{p-1}b_i+1<\sum\limits_{i=0}^{p}a_i+2k\right\}.
\]
The inner product {\em(\ref{InnerDiskSeries})} can be then written as
\[
[f,g]\!=\!\frac{\theta^{1-l}}{\pi}\Gamma\!\!\left[\!\begin{array}{l}(a_p)\\(b_{p-1})\end{array}\!\right]
\int\limits_{\Delta_{\sqrt{\theta}}}\!\!\!f^{(l)}(z)\overline{g^{(l)}(z)}
\G^{p+1,0}_{p+1,p+1}\left(\frac{|z|^2}{\theta}\vline\begin{array}{l}(a_p)+l-1,l\\
(b_{p-1})+l-1,0,0\end{array}\!\!\!\right)\!d\sigma+
\]
\begin{equation}\label{InnerDisk1}
+\sum\limits_{k=0}^{l-1}\frac{(b_1)_k\ldots(b_{p-1})_k
\theta^k}{(a_1)_k\ldots(a_p)_k k!}f^{(k)}(0)\overline{g^{(k)}(0)},
\end{equation}
with $l=\max(l_1,l_2)$.

If $l_1=0$, then the inner product {\em(\ref{InnerDiskSeries})} can
also be written as
\begin{equation}\label{InnerDisk2}
[f,g]\!=\!\frac{1}{\pi}\Gamma\!\!\left[\!\begin{array}{l}(a_p)\\(b_{p-1})\end{array}\!\right]
\int\limits_{\Delta_{\sqrt{\theta}}}\!\!\!\left((Dz)^lf(z)\right)\overline{\left((Dz)^lg(z)\right)}
\frac{1}{|z|^{2}}
\G^{p+2l,0}_{p+2l,p+2l}\left(\frac{|z|^2}{\theta}\vline\begin{array}{l}(a_p),(2)\\
(b_{p-1}),(1)\end{array}\!\!\!\right)\!d\sigma.
\end{equation}
In particular if $l_1=l_2=0$
\begin{equation}\label{InnerDisk}
[f,g]\!=\!\frac{1}{\theta\pi}\Gamma\!\!\left[\!\begin{array}{l}(a_p)\\(b_{p-1})\end{array}\!\right]
\int\limits_{\Delta_{\sqrt{\theta}}}\!\!\!f(z)\overline{g(z)}
\G^{p,0}_{p,p}\left(\frac{|z|^2}{\theta}\vline\begin{array}{l}(a_p)-1\\0,(b_{p-1})-1\end{array}
\!\!\!\right)\!d\sigma.
\end{equation}
\end{theo}

The spaces $\B^{\,\theta}_{(a_p),(b_{p-1})}$ contain, as particular
cases, the classical Bergman ($q=0$, $a=2$), the Hardy $\H_2$ ($q=0$,
$a=1$) and the Bergman-Selberg ($q=0$, $a>1$) spaces, as well as many
others (see \cite{Karp4} for other examples).

\dept{3. Topologies of hypergeometric kernel spaces}

In this section we consider a generalization of the kernel
(\ref{pFq-def}) given by
\begin{equation}\label{Hyperkern}
K(z,\overline{u})=\sum\limits_{k=0}^{\infty}
\frac{\Gamma(B_1k+b_1)\cdots\Gamma(B_qk+b_q)}{\Gamma(A_1k+a_1)\cdots\Gamma(A_pk+a_p)}
(z\overline{u})^k.
\end{equation}
The series on the right is called the Wright psi-function
\cite{Wright1,Wright2}. For simplicity we limit ourselves to the
Hilbert space case $A_i,B_i,a_i,b_i>0$. Formula (\ref{GMellin}) is
generalized as follows (cf.\cite{Mathai,Prud3})
\begin{equation}\label{HMellin}
\int\limits_{0}^{\infty}\!\!x^{s-1}
{\mathrm{H}}^{q,0}_{p,q}\left(x~\vline\begin{array}{l}(a_p,A_p)\\(b_q,B_q)\end{array}\right)dx=
\frac{\Gamma(b_1+B_1s)\ldots\Gamma(b_q+B_qs)}{\Gamma(a_1+A_1s)\ldots\Gamma(a_p+A_ps)}.
\end{equation}
Here ${\mathrm{H}}^{q,0}_{p,q}$ is Fox's H-function
\cite{Braaksma,Mathai,Prud3}.  Thus one can construct the analogues
of Theorems~\ref{KreinPlane} and \ref{KreinDisk} for the spaces
generated by the kernel (\ref{Hyperkern}).  The purpose of this
section, however, is different.  We show that while the geometries of
the spaces generated by the kernel (\ref{Hyperkern}) depend on all
parameters $A_i, B_i, a_i, b_i$ their topologies are only dependent
on three numbers
\[
\alpha=\sum\limits_{i=1}^{q}b_i-\sum\limits_{i=1}^{p}a_i+\frac{p-q}{2},
\]
\[
\mu=\sum\limits_{i=1}^{q}B_i-\sum\limits_{i=1}^{p}A_i\geq 0
\]
and
\[
\nu=\prod\limits_{i=1}^{p}A_i^{-A_i}\prod\limits_{i=1}^{q}B_i^{B_i}>0.
\]
Indeed, the norms in the Hilbert spaces $H_{K^{(\varepsilon)}}$,
$\varepsilon\in\{1,2\}$ induced by the kernels $K^{(1)}$, $K^{(2)}$
of the form (\ref{Hyperkern}) are equivalent when
\[
\frac{\Gamma(B_1^{(1)}k+b_1^{(1)})\cdots\Gamma(B_q^{(1)}k+b_q^{(1)})}
{\Gamma(A_1^{(1)}k+a_1^{(1)})\cdots\Gamma(A_p^{(1)}k+a_p^{(1)})}\sim
\frac{\Gamma(B_1^{(2)}k+b_1^{(2)})\cdots\Gamma(B_q^{(2)}k+b_q^{(2)})}
{\Gamma(A_1^{(2)}k+a_1^{(2)})\cdots\Gamma(A_p^{(2)}k+a_p^{(2)})},~~k\to\infty.
\]
From Stirling's formula one easily derives
\begin{equation}\label{Gamfracas}
\frac{\Gamma(B_1k+b_1)\cdots\Gamma(B_qk+b_q)}{\Gamma(A_1k+a_1)\cdots\Gamma(A_pk+a_p)}
\sim k^{\alpha}e^{-\mu{k}}k^{\mu{k}}\nu^k,~~k\to\infty,
\end{equation}
which proves our conjecture about the topologies of $H_K$. Thus we
obtain two ''model'' spaces for two distinct cases $\mu>0$ and
$\mu=0$.  If $\mu>0$ the space induced by (\ref{Hyperkern}) is
topologically equivalent to the space of entire functions equipped
with the norm
\[
\int\limits_{\C}|f^{(l)}(z)|^2
\exp\left(-\mu\nu^{-\mu}|z|^{2/\mu}\right)|z|^{(2\alpha+4l+2\mu{l}+1)/\mu-2}d\sigma+
\sum\limits_{k=0}^{l-1}|f^{(k)}(0)|^2,
\]
where $l=\max(0,[(-\alpha-\frac{1}{2})/\mu]+1)$. In particular, when
$\alpha>-1/2$ the norm reduces to
\[
\frac{2\alpha+1}{2\pi\mu}\int\limits_{\C}|f(z)|^2
\exp\left(-\mu\nu^{-\mu}|z|^{2/\mu}\right)|z|^{(2\alpha+1)/\mu-2}d\sigma,
\]
and the reproducing kernel is
$E_{\mu,\alpha+1/2}(\nu^{-1}\mu^{\mu}z\overline{u})$, where
\[
E_{\mu,\alpha}(t)=\sum\limits_{k=0}^{\infty}\frac{t^k}{\Gamma(\mu{k}+\alpha)}
\]
is the two-parameter Mittag-Leffler function \cite{DJ}.

If $\mu=0$ the space induced by (\ref{Hyperkern}) is topologically
equivalent to the space of functions analytic inside
$\Delta_{\sqrt{\nu}}$ equipped with the norm
\[
\int\limits_{\Delta_{\sqrt{\nu}}}|f^{(l)}(z)|^2
(1-|z|^2/\nu)^{2l-\alpha-1}d\sigma+
\sum\limits_{k=0}^{l-1}|f^{(k)}(0)|^2,
\]
where $l=\max(0,[\alpha]+1)$.  In particular, when $\alpha<0$ it is
the Bergman-Selberg (or weighted Bergman) space with the norm
\[
\frac{-\alpha}{\pi\nu}\int\limits_{\Delta_{\sqrt{\nu}}}|f(z)|^2(1-|z|^2/\nu)^{-\alpha-1}d\sigma
\]
and reproducing kernel
\[
K(z,\overline{u})=\left(1-\frac{z\overline{u}}{\nu}\right)^{\alpha-1}.
\]

\dept{4. Analytic continuation from the positive half-line}

In this section we use a hypergeometric kernel space to derive an
analytic extension formula for functions defined on $\R^+$.  As a
by-product several new area integrals involving Bessel functions are
explicitly evaluated.

We shall use the theory of best approximation in reproducing kernel
Hilbert spaces (RKHS) developed by S. Saitoh as set forth in
\cite{SaitohBook}. Let $E$ be an arbitrary set, and let $H_K$ be a
RKHS comprising functions $E\to\C$ and admitting the reproducing
kernel $K$. For $X\subset{E}$ we consider a Hilbert space $H(X)$
comprising functions $X\to\C$.  We additionally assume that
restrictions $f|_X$ of $f\in{H_K}$ belong to the Hilbert space $H(X)$
and the restriction operator $T$ is continuous from $H_K$ into
$H(X)$.  In this setting Saitoh has proved the following result.
\begin{theo}\label{BestApproxTheo}
For a given $d\in{H(X)}$ there exists $\tilde{f}\in{H_K}$ such that
\begin{equation}\label{BestApprox}
\inf\limits_{f\in{H_K}}\|Tf-d\|_{H(X)}=\|T\tilde{f}-d\|_{H(X)}
\end{equation}
if and only if $T^*d\in{H_k}$, where RKHS $H_k$ is induced by the
kernel
\begin{equation}\label{BestKernel}
k(p,q)=\left(T^*TK(\cdot,q),T^*TK(\cdot,p)\right)_{H_K}~~~\mathrm{on}~~
E\times{E}.
\end{equation}
Among all solutions of {\em(\ref{BestApprox})} there exists a unique
extremal $f^*$ with the minimum norm in $H_K$ which is given by
\begin{equation}\label{Extremal}
f^*(p)=(T^*d,T^*TK(\cdot,p))_{H_k}.
\end{equation}
\end{theo}

This theorem has been applied in \cite{BS} to derive an analytic
extension formula for function defined on $\R$.  In \cite{Karp5} the
author found some analytic continuability criteria for functions on
$\R$, $\R^+$, and $[-1,1]$ without the application of
Theorem~\ref{BestApproxTheo}.

To obtain a result for functions on $\R^+$ similar to that of
\cite{BS} consider the space $\H^{\alpha,\nu}_K$ induced by the
kernel
\begin{equation}\label{Ikernel}
K_{\alpha,\nu}(z,\overline{u})=\frac{1}{\Gamma(\nu+1)}{_{0}\mathrm{F}_{1}}
\left(-;\nu+1;\frac{\alpha^{2}z\overline{u}}{4}\right)=
\sum\limits_{k=0}^{\infty}\frac{\alpha^{2k}(z\overline{u})^k}{4^k\Gamma(k+\nu+1)k!}=
\left(\frac{2}{\alpha}\right)^{\nu}(z\overline{u})^{-\nu/2}\I_{\nu}(\alpha\sqrt{z\overline{u}}),
\end{equation}
where $\I_\nu$ is the modified Bessel function (cf.
\cite{Bat2,Prud2}), $\nu>-1$ and $\alpha>0$. By
Theorem~\ref{KreinPlane} this space consists of entire functions with
finite norms
\begin{equation}\label{Knorm}
\|f\|^2_{\nu,\alpha}=\frac{\alpha^{\nu+2}}{2^{\nu+1}\pi}\int\limits_{\C}|f(z)|^2|z|^{\nu}K_{\nu}(\alpha|z|)d\sigma,
\end{equation}
where $K_\nu$ is the MacDonald function (cf. \cite{Bat2,Prud2}). We
consider the restriction operator $T: \H^{\alpha,\nu}_K\to
L_2(\R^+;W_{\alpha,\nu})$, $W_{\alpha,\nu}(t)=e^{-\alpha{t}}t^{\nu}$.
\begin{lemma}
The restriction operator $T$ is bounded.
\end{lemma}
\textbf{Proof.} From the formula
\begin{equation}\label{Key}
\int\limits_{0}^{\infty}J_{\nu}(zt)J_{\nu}(ut)e^{-{\gamma}t^2}tdt=
\frac{1}{2\gamma}\exp\left[-\frac{z^2+u^2}{4\gamma}\right]\I_{\nu}\left(\frac{zu}{2\gamma}\right),
~~\gamma>0,
\end{equation}
found for example in \cite{Bat2,Prud2}, we get the following
representation for the kernel
\begin{equation}\label{Krepresent}
K_{\alpha,\nu}(z,\overline{u})=\left(\frac{2}{\alpha}\right)^{\nu}
e^{{\alpha}z/2}e^{{\alpha}\overline{u}/2}\int\limits_{0}^{\infty}(zt)^{-\nu/2}J_\nu(\sqrt{zt})
(\overline{u}t)^{-\nu/2}J_\nu(\sqrt{\overline{u}t})e^{-t/(2\alpha)}t^{\nu}dt.
\end{equation}
In these formulae $J_\nu$ is the Bessel function.  According to the
theory developed in \cite{SaitohBook}, this representation of the
kernel implies that each $f\in\H^{\alpha,\nu}_K$ is expressible in
the form
\begin{equation}\label{frepresent}
f(z)=\left(\frac{2}{\alpha}\right)^{\nu} e^{{\alpha}z/2}
\int\limits_{0}^{\infty}(zt)^{-\nu/2}J_\nu(\sqrt{zt})F(t)e^{-t/(2\alpha)}t^{\nu}dt
\end{equation}
for some uniquely determined function $F$ satisfying
\[
\int\limits_{0}^{\infty}|F(t)|^2e^{-t/(2\alpha)}t^{\nu}dt<\infty,
\]
and we have the isometric identity
\[
\|f\|^2_{\nu,\alpha}=\left(\frac{2}{\alpha}\right)^{\nu}\int\limits_{0}^{\infty}|F(t)|^2e^{-t/(2\alpha)}t^{\nu}dt.
\]
On the other hand from the unitarity of the Hankel transform in
$L_2(\R^+)$ it follows by a simple change of variable that the
operator
\[
A:~
g\to\frac{1}{2}\int\limits_{0}^{\infty}J_{\nu}\left(\sqrt{\cdot~t}\right)g(t)dt
\]
is also unitary in $L_2(\R^+)$.  Hence from (\ref{frepresent})
\[
\int\limits_{0}^{\infty}|f(t)|^2W_{\alpha,\nu}(t)dt=\frac{2^{2\nu+2}}{\alpha^{2\nu}}
\int\limits_{0}^{\infty}\left|F(t)e^{-t/(2\alpha)}t^{\nu/2}\right|^2dt
\]\[
\leq4\left(\frac{2}{\alpha}\right)^{2\nu}\int\limits_{0}^{\infty}|F(t)|^2e^{-t/(2\alpha)}t^{\nu}dt=
4\left(\frac{2}{\alpha}\right)^{\nu}\|f\|^2_{\nu,\alpha}.~~\square
\]
Now we consider the existence of the best approximation $f^*$
\[
\|Tf^*-F\|_{L_2(\R^+,W_{\alpha,\nu})}=\inf\limits_{f\in\H^{\alpha,\nu}_K}\|Tf-F\|_{L_2(\R^{+},W_{\alpha,\nu})}=0.
\]
The last equality is due to the density in $L_2(\R^+;W_{\alpha,\nu})$
of restrictions of functions from $\H^{\alpha,\nu}_K$ to $\R^{+}$.
$T^*$ is expressible in the form (by definition of the adjoint
operator):
\[
(T^*F(\cdot))(u)=((T^*F)(\cdot),K_{\alpha,\nu}(\cdot,\overline{u}))_{\H^{\alpha,\nu}_K}
=(F(\cdot),TK_{\alpha,\nu}(\cdot,\overline{u}))_{L_2(\R^+;W_{\alpha,\nu})}=
\]\[
=\int\limits_{0}^{\infty}\!\!F(t)\overline{K_{\alpha,\nu}(t,\overline{u})}W_{\alpha,\nu}(t)dt=
\int\limits_{0}^{\infty}\!\!F(t)K_{\alpha,\nu}(t,u)W_{\alpha,\nu}(t)dt.
\]
Then for $s\in\R^+$
\[
\left(TT^*TK_{\alpha,\nu}(\cdot,\overline{z})\right)(s)=
\int\limits_{0}^{\infty}\!\!K_{\alpha,\nu}(t,\overline{z})K_{\alpha,\nu}(t,s)W_{\alpha,\nu}(t)dt
\]\[
=\left(\frac{2}{\alpha}\right)^{2\nu}(s\overline{z})^{-\nu/2}
\int\limits_{0}^{\infty}\!\I_{\nu}(\alpha\sqrt{t\overline{z}})
\I_{\nu}(\alpha\sqrt{ts})e^{-{\alpha}t}dt
\]\[
=
\frac{1}{\alpha}\left(\frac{2}{\alpha}\right)^{2\nu}(s\overline{z})^{-\nu/2}
e^{\alpha(\overline{z}+s)/4}
\I_{\nu}\left(\frac{\alpha}{2}\sqrt{s\overline{z}}\right)=
\alpha^{-\nu-1}e^{\alpha(\overline{z}+s)/4}K_{\alpha/2,\nu}(s,\overline{z}).
\]
Hence for the kernel $k$ from Theorem~\ref{BestApproxTheo} we
calculate
\[
k_{\alpha,\nu}(z,\overline{u})=
\left(T^*TK_{\alpha,\nu}(\cdot,\overline{u}),T^*TK_{\alpha,\nu}(\cdot,\overline{z})\right)_{\H^{\alpha,\nu}_K}=
\left(TK_{\alpha,\nu}(\cdot,\overline{u}),TT^*TK_{\alpha,\nu}(\cdot,\overline{z})\right)_{L_2(\R^+,W_{\alpha,\nu})}
\]\[
=
\alpha^{-\nu-1}e^{\alpha{z}/4}\int\limits_{0}^{\infty}\!\!K_{\alpha,\nu}(s,\overline{u})
K_{\alpha/2,\nu}(s,z)e^{\alpha{s}/4}W_{\alpha,\nu}(s)ds
\]\[
=\frac{1}{\alpha}\left(\frac{2}{\alpha}\right)^{3\nu}
e^{\alpha{z}/4}(z\overline{u})^{-\nu/2}
\int\limits_{0}^{\infty}\!\!\I_{\nu}(\alpha\sqrt{s\overline{u}})
\I_{\nu}\left(\frac{\alpha}{2}\sqrt{sz}\right)e^{-3\alpha{s}/4}ds
\]\[
=\frac{4}{3\alpha^2}\left(\frac{2}{\alpha}\right)^{3\nu}
e^{\alpha{z}/3}e^{\alpha\overline{u}/3}(z\overline{u})^{-\nu/2}
\I_{\nu}\left(\frac{\alpha\sqrt{z\overline{u}}}{3}\right)=
\left(\frac{4}{3}\right)^{\nu+1}\alpha^{-2\nu-2}
e^{\alpha{z}/3}e^{\alpha\overline{u}/3}K_{\alpha/3,\nu}(z,\overline{u}).
\]
From the general theory of reproducing kernels
\cite{Ar2,SaitohBook} and (\ref{Knorm}) we see that the space
$\H^{\alpha,\nu}_k$ induced by the kernel $k_{\alpha,\nu}$
comprises entire functions with finite norms
\[
\frac{\alpha^{3\nu+4}}{3\pi2^{3\nu+3}}
\int\limits_{\C}|f(z)|^2\exp\left(-\frac{2}{3}\alpha\Re{z}\right)|z|^{\nu}K_{\nu}\left(\frac{\alpha|z|}{3}\right)d\sigma.
\]
This brings us to the main result of this section.  In the
following theorem analytic continuation is understood in the sense
that the extended and the original functions may differ on a set
of zero measure.
\begin{theo}\label{ContinTheo}
A function $F\in L_2(\R^+;W_{\alpha,\nu})$ can be analytically
continued to the entire complex plane to a function $f$ from the
space $\H^{\alpha,\nu}_K$ if and only if
\begin{equation}\label{Criterion}
\int\limits_{\C}\left|\int\limits_{0}^{\infty}\!\!F(t)
\I_{\nu}(\alpha\sqrt{tz})e^{-\alpha{t}}t^{\nu/2}dt\right|^2
e^{-\frac{2}{3}\alpha\Re{z}}K_{\nu}\left(\frac{\alpha|z|}{3}\right)d\sigma<\infty.
\end{equation}
The analytic continuation is given by
\begin{equation}\label{Contin}
f(z)=\frac{\alpha^{3}}{3\pi2^{3}}e^{\alpha{z}/4}z^{-\nu/2}
\int\limits_{\C}\I_{\nu}\left(\frac{\alpha}{2}\sqrt{z\overline{u}}\right)
e^{-\frac{\alpha}{4}{u}-\frac{\alpha}{6}\Re{u}}K_{\nu}\left(\frac{\alpha|u|}{3}\right)
\int\limits_{0}^{\infty}\!\!
F(t)\I_{\nu}(\alpha\sqrt{tu})e^{-\alpha{t}}t^{\nu/2}dtd\sigma(u).
\end{equation}
\end{theo}

Formula (\ref{Contin}) contains many interesting special cases.
Taking $F(t)=1$ we get
\[
e^{-\alpha{z}/4}z^{\nu/2}=\frac{\alpha^{2}}{3\pi2^{3+\nu}}
\int\limits_{\C}u^{\frac{\nu}{2}}\I_{\nu}\left(\frac{\alpha}{2}\sqrt{z\overline{u}}\right)
e^{-\frac{1}{6}\alpha\Re{u}}K_{\nu}\left(\frac{\alpha|u|}{3}\right)d\sigma(u).
\]
In particular, when $\alpha=2$ and $\nu=\frac{1}{2}$ the last
equality reduces to
\[
z^{\frac{1}{2}}e^{-\frac{z}{2}}=\frac{1}{4\pi\sqrt{3}}\int\limits_{\C}
u^{\frac{1}{2}}\sinh\left(\sqrt{z\overline{u}}\right)|u|^{-1}e^{-\frac{2}{3}|u|-\frac{1}{3}\Re{u}}d\sigma(u),
\]
which can be verified directly. If we put $F(t)=t^n$, $n\in\N$,
then evaluation of the integral over $\R^+$ in (\ref{Contin})
leads to
\begin{equation}\label{Lag2power}
z^{n+\nu/2}e^{-\alpha{z}/4}=\frac{n!}{3\pi2^{\nu+3}\alpha^{n-2}}
\int\limits_{\C}u^{\frac{\nu}{2}}\I_{\nu}\left(\frac{\alpha}{2}\sqrt{z\overline{u}}\right)
\L^{\!\nu}_n\left(-\frac{\alpha{u}}{4}\right)
e^{-\frac{1}{6}\alpha\Re{u}}K_{\nu}\left(\frac{\alpha|u|}{3}\right)
d\sigma(u),
\end{equation}
where $\L^{\!\nu}_n$ is the Laguerre polynomial \cite{Bat2,Prud2}.
We can derive a counterpart of (\ref{Lag2power}) by using (cf.
\cite{Prud2})
\[
u^n=\frac{4^{n}n!(\nu+1)_n}{\alpha^n}\sum\limits_{k=0}^{n}
\frac{(-1)^{n-k}}{(n-k)!(\nu+1)_k}\L^{\!\nu}_k\left(-\frac{\alpha{u}}{4}\right).
\]
Substituting this identity into (\ref{Lag2power}) yields
\begin{equation}\label{Power2Lag}
\int\limits_{\C}u^{n+\frac{\nu}{2}}\I_{\nu}\left(\frac{\alpha}{2}\sqrt{z\overline{u}}\right)
e^{-\frac{1}{6}\alpha\Re{u}}K_{\nu}\left(\frac{\alpha|u|}{3}\right)
d\sigma(u)=\frac{(-1)^{n}3\pi2^{\nu+3}4^{n}n!}{\alpha^{n+2}}z^{\nu/2}e^{-\alpha{z}/4}
\L^{\!\nu}_n(\alpha{z}).
\end{equation}
One more simple case is $F(t)=t^{-\nu/2}\I_{\nu}(\sqrt{t})$, which
results in
\[
e^{-\alpha{z}/4}\I_{\nu}(\sqrt{z})=\frac{\alpha^{2}e^{1/(4\alpha)}}{3\pi2^{3}}
\int\limits_{\C}\I_{\nu}\left(\frac{\alpha}{2}\sqrt{z\overline{u}}\right)
\I_{\nu}\left(\frac{\sqrt{u}}{2}\right)
e^{-\frac{1}{6}\alpha\Re{u}}K_{\nu}\left(\frac{\alpha|u|}{3}\right)
d\sigma(u).
\]

Further following \cite{BS}, for {\em any} $F\in
L_2(\R^+;W_{\alpha,\nu})$ we can construct a sequence
$\{f_n\}\in\H^{\alpha,\nu}_K$ such that
\[
\lim\limits_{n\rightarrow\infty}\|Tf_n-F\|_{L_2(\R^+;W_{\alpha,\nu})}=0.
\]
The general results for integral transforms set forth in
\cite{SaitohBook} imply that the images $f$ of the transform
\[
(T^*F)(z)=(F(\cdot),TK_{\alpha,\nu}(\cdot,\overline{z}))_{L_2(\R^+;W_{\alpha,\nu})}=
\int\limits_{0}^{\infty}\!\!F(t)K_{\alpha,\nu}(t,z)W_{\alpha,\nu}(t)dt
\]
for $F\in L_2(\R^+;W_{\alpha,\nu})$  form the RKHS with reproducing
kernel
\begin{equation}\label{rformula}
r(z,\overline{u})=
\left(TK_{\alpha,\nu}(\cdot,\overline{u}),TK_{\alpha,\nu}(\cdot,\overline{z})\right)_{L_2(\R^+;W_{\alpha,\nu})}=
\alpha^{-\nu-1}e^{\alpha{z}/4}e^{\alpha\overline{u}/4}K_{\alpha/2,\nu}(z,\overline{u})
\end{equation}
as has been previously computed. The norm in $\H_{r}$ equals then
\[
\|f\|^2_{\H_{r}}=\frac{\alpha^{2\nu+3}}{2^{2\nu+3}\pi}
\int\limits_{\C}|f(z)|^2e^{-\frac{\alpha}{2}\Re{z}}|z|^{\nu}K_{\nu}(\alpha|z|/2)d\sigma.
\]
From representation of the kernel (\ref{rformula}) we infer that
the functions
\begin{equation}\label{phidefin}
\varphi_{k,\nu}(z)=\frac{\alpha^{k}e^{\alpha{z}/4}z^k}{4^k(\alpha^{\nu+1}\Gamma(k+\nu+1)k!)^{1/2}}
\end{equation}
form a complete orthonormal system in $\H_r$.  Note also that
$\varphi_{k,\nu}\in\H^{\alpha,0}_k$ for any $k\in\N_0$. Consider
the partial Fourier sum
\[
\tilde{f}_N(z)=\sum\limits_{k=0}^{N}a_k\varphi_{k,\nu}(z),~~~a_k=(T^*F,\varphi_{k,\nu})_{\H_r}.
\]
$\{\tilde{f}_N\}$ converges to $T^*F$ in $\H_r$ as
$N\rightarrow\infty$. We can construct the sequence $f^*_N$ such
that $T^*Tf^*_N=\tilde{f}_N$ by Theorem~\ref{ContinTheo}.  First
we need $\varphi^*_{k,\nu}$ such that
$T^*T\varphi^*_{k,\nu}=\varphi_{k,\nu}$. From (\ref{Contin})
\[
\varphi^*_{k,\nu}(z)=
\frac{\alpha^{k+3}e^{\alpha{z}/4}z^{-\nu/2}}{3\pi2^{2k+3}(\alpha^{\nu+1}\Gamma(k+\nu+1)k!)^{1/2}}
\int\limits_{\C}\I_{0}\left(\frac{\alpha}{2}\sqrt{z\overline{u}}\right)
u^ke^{-\frac{\alpha}{6}\Re{u}}K_{0}\left(\frac{\alpha|u|}{3}\right)
d\sigma(u)
\]\[
=\frac{\alpha^{\frac{1}{2}}(\alpha{z})^{-\nu/2}(-1)^{k}(k!)^{1/2}}{\Gamma(k+\nu+1)^{1/2}}
\L^{\!0}_k(\alpha{z}),
\]
where (\ref{Power2Lag}) has been utilized with $\nu=0$.  Now
compute $a_k$:
\[
a_k=\frac{\alpha^{k+\nu/2+5/2}}{2^{\nu+3}\pi4^k(\Gamma(k+\nu+1)k!)^{1/2}}
\int\limits_{0}^{\infty}\!\!F(t)e^{-\alpha{t}}t^{\nu/2}A^{\alpha}_{\nu}(t)dt,
\]
\[
A^{\alpha}_{\nu}(t)=\int\limits_{\C}e^{-\alpha{u}/4}\overline{u}^k
u^{-\nu/2}\I_{\nu}(\alpha\sqrt{tu})
|u|^{\nu}K_{\nu}(\alpha|u|/2)d\sigma.
\]
The last integral cannot be evaluated using (\ref{Contin}).  Finally
we have
\[
f^*_N=\sum\limits_{k=0}^{N}a_k\varphi^*_{k,\nu},
\]
and since $T^*$ is an isometry from $L_2(\R^+,W_{\alpha,\nu})$ onto
$\H_r$, we obtain
\[
\lim\limits_{N\rightarrow\infty}\|f^*_N-F\|_{L_2(\R^+,W_{\alpha,\nu})}=0.
\]
It follows from the results of \cite{Karp4,Karp5}  that the
functions $\varphi^*_{k,\nu}$ satisfy the orthogonality relation
given by
\[
\int\limits_{\C}\!\!\ \varphi^*_{k,\nu}(z) \overline{\varphi^*_{m,\nu}(z)}
\exp\!\left(\frac{2\alpha\Re{z}}{\theta-1}\right)
|z|^{\nu}K_{0}\left(\frac{2\alpha\sqrt{\theta}|z|}{\theta-1}\right)d\sigma=\delta_{k,m}
\frac{\pi(\theta-1)\theta^k(-1)^{m+k}(m!k!)^{1/2}}
{2\alpha^{\nu+1}(\Gamma(m+\nu+1)\Gamma(k+\nu+1))^{1/2}}
\]
with arbitrary $\theta>1$.  Note that the functions $\varphi^*_k$
from \cite{BS} can be found explicitly in terms of the Hermite
polynomials.  They satisfy an orthogonality relation easily
derivable from the orthogonality relation for the Hermite
polynomials given in \cite{Karp5}.

\dept{Acknowledgements.} This research is carried out with financial
support from Russian Basic Research Foundation (grant no.
02-01-000-28) and the "Russian Universities" program (grant no.
UR.04.01.016).

\end{document}